\documentclass[fleqn,11pt,twoside]{article}

\usepackage{amsthm,amsthm,amssymb, color, xcolor,epsfig, graphics, subfigure}

\usepackage{amsmath, graphicx, latexsym, lscape}

\usepackage{mathrsfs}
\usepackage{bbm}
\usepackage{physics}
\numberwithin{equation}{section}

\makeatletter
\newcommand{\copyrightnote}[2]{{\renewcommand{\thefootnote}{}
 \footnotetext{\small\it
\begin{flushleft}
 \copyright \ #1   #2  
\end{flushleft}}}}

\newcommand{\Name}[1]{\begin{flushleft}
                       \LARGE \bf #1
                       \end{flushleft}\vspace{-3mm}}

\newcommand{\Author}[1]{\begin{flushleft}
                       \it #1 \end{flushleft}}

\newcommand{\Address}[1]{\begin{flushleft}
                       \it #1 \end{flushleft}}

\newcommand{\Date}[1]{\begin{flushleft}
                      \small  \it #1 \end{flushleft}}

\newcommand{\evenhead}{Author \ name}
\newcommand{\oddhead}{Article \ name}

\renewcommand{\@evenhead}{
\hspace*{-3pt}\raisebox{-15pt}[\headheight][0pt]{\vbox{\hbox to \textwidth
{\thepage \hfil \evenhead}\vskip4pt \hrule}}}
\renewcommand{\@oddhead}{
\hspace*{-3pt}\raisebox{-15pt}[\headheight][0pt]{\vbox{\hbox to \textwidth
{\oddhead \hfil \thepage}\vskip4pt\hrule}}}
\renewcommand{\@evenfoot}{}
\renewcommand{\@oddfoot}{}

\long\def\@makecaption#1#2{%
  \vskip\abovecaptionskip
  \sbox\@tempboxa{\small \textbf{#1.}\ \ #2}%
  \ifdim \wd\@tempboxa >\hsize
    {\small \textbf{#1.}\ \ #2}\par
  \else
    \global \@minipagefalse
    \hb@xt@\hsize{\hfil\box\@tempboxa\hfil}%
  \fi
  \vskip\belowcaptionskip}

\newcommand{\JNMPnumberwithin}[3][\arabic]{%
  \@ifundefined{c@#2}{\@nocounterr{#2}}{%
    \@ifundefined{c@#3}{\@nocnterr{#3}}{%
      \@addtoreset{#2}{#3}%
      \@xp\xdef\csname the#2\endcsname{%
        \@xp\@nx\csname the#3\endcsname .\@nx#1{#2}}}}%
}

\newcommand{\resetfootnoterule} {
  \renewcommand\footnoterule{%
  \kern-3\p@
  \hrule\@width.4\columnwidth
  \kern2.6\p@}
}

\renewcommand{\footnoterule}{}

\makeatother

\theoremstyle{definition}

\begin{document}

\renewcommand{\evenhead}{ {\LARGE\textcolor{blue!10!black!40!green}{{\sf \ \ \ ]ocnmp[}}}\strut\hfill A M Grundland}
\renewcommand{\oddhead}{ {\LARGE\textcolor{blue!10!black!40!green}{{\sf ]ocnmp[}}}\ \ \ \ \  On $k$-waves solutions of quasilinear systems of PDEs }

\thispagestyle{empty}
\newcommand{\FistPageHead}[3]{
\begin{flushleft}
\raisebox{8mm}[0pt][0pt]
{\footnotesize \sf
\parbox{150mm}{{Open Communications in Nonlinear Mathematical Physics}\ \ \ \ {\LARGE\textcolor{blue!10!black!40!green}{]ocnmp[}}
\quad Special Issue 1, 2024\ \  pp
#2\hfill {\sc #3}}}\vspace{-13mm}
\end{flushleft}}

\FistPageHead{1}{\pageref{firstpage}--\pageref{lastpage}}{ \ \ }

\strut\hfill

\strut\hfill

\copyrightnote{The author(s). Distributed under a Creative Commons Attribution 4.0 International License}

\begin{center}
{  {\bf This article is part of an OCNMP Special Issue\\ 
\smallskip
in Memory of Professor Decio Levi}}
\end{center}

\smallskip

\Name{On $k$-wave solutions of quasilinear systems of partial differential equations}

\Author{A.M. Grundland$^{\,1}$,$^{\,2}$}

\Address{$^{1}$ Centre de Recherches Mathématiques, Université de Montréal, \\Succ. Centre-Ville, CP 6128, Montréal (QC) H3C 3J7, Canada\\[2mm]
$^{2}$ Département de Mathématiques et d'Informatique, \\Université du Québec, CP 500, Trois-Rivières (QC)\\ G9A 5H7, Canada}

\Date{Received May 21, 2023; Accepted October 6, 2023}

\setcounter{equation}{0}

\begin{abstract}

\noindent 
In this paper, we establish a relation between two seemingly unrelated concepts for solving first-order hyperbolic quasilinear systems of partial differential equations in many dimensions. These concepts are based on a variant of the conditional symmetry method and on the generalized method of characteristics. We present the outline of recent results on multiple Riemann wave solutions of these systems. An auxiliary result concerning a modification of the Frobenius theorem for integration is used. We apply this result in order to show that the conditional symmetry method can deliver larger classes of multiple Riemann wave solutions, through a simpler procedure, than the one obtained from the generalized method of characteristics. We demonstrate that solutions can be interpreted physically as a superposition of $k$ single waves. These theoretical considerations are illustrated by examples of hydrodynamic-type systems in $(n+1)$-dimensions.

\end{abstract}

\label{firstpage}


\section{Introduction}

This work is dedicated to the memory of professor Decio Levi (University of Roma Tre), a friend and collaborator whose support and knowledgeable and generous advice I could always count on during the many years of our aquaintance. The subject of the present article was inspired by Decio Levi's encouragement, which resulted in a series of lectures that I gave at the University of Roma Tre in the spring of 2016.
	
	Over the last two centuries, a number of authors (e.g. \cite{Lie,OlverLie,Stephani,Winternitz,Ibragimov} and references therein) have undertaken a systematic study of the group invariant solutions of nonlinear partial differential equations (PDEs). This study has often required a detailed analysis of the subgroup structures of a given Lie group which arises as a symmetry group of some system of PDEs. A number of attempts to generalize this subject and to develop its applications can be found in the literature (see e.g. \cite{Ibragimov,OlverSymm,BlumanCole} and references therein). Of particular interest from a physical point of view has been the development of the theory of conditional symmetries which is involved in the process of extending the classical Lie theory of symmetries of PDEs. This approach consists essentially of augmenting the initial system of PDEs with certain first-order differential constraints (DCs) for which a symmetry criterion is applied. This leads to an overdetermined system of PDEs that can admit, in some cases, a larger family of Lie point symmetries than the initial system. Such Lie point symmetries can be used to find and construct particular types of solutions of the original system of PDEs \cite{OlverSymm,BlumanCole,GrundlandRideau,GrundlandMartinaRideau,OlverRosenan,OlverVorobiev,Zidowitz}. If, in an attempt to obtain larger classes of solutions, one adds DCs giving rise to a first-order system of PDEs in a form which is compatible (in a sense which will be explained in the next sections) with the initial system of PDEs, then one obtains the conditional symmetry formalism \cite{OlverVorobiev,Zidowitz}. In this formalism, the DCs are described through a Lie algebra of vector fields, whose elements are Lie point symmetries of the initial system of PDEs leading to solutions satisfying the imposed DCs. To obtain an easily applicable formalism, Lie algebras of conditional symmetries are assumed to obey some assumptions, e.g. one considers Abelian Lie algebras admitting a basis of a particular type \cite{GrundlandRideau,GrundlandMartinaRideau}.
	
	In this paper, we focus on the investigation and construction of multiple Riemann wave solutions obtained via the conditional symmetry method and a comparison of these results with the ones obtained through the generalized method of characteristics established in \cite{Burnat,PeradzynskiNonlin,GrundlandInvolutivity}.
	
	The method of characteristics is used to find a certain class of solutions of nonlinear hyperbolic systems of PDEs and to investigate the existence and construction of the Riemann wave and Riemann $k$-wave solutions. The single wave solution was first introduced by S. Poisson in 1808 in connection with PDEs describing an ideal compressible isothermal fluid flow. He constructed a solution assuming that it had an implicit form between the dependent and independent $x$ variables with the freedom of one arbitrary function $F$, i.e. $u=F\left((u+u_0)t-x\right)$. This idea was substantially generated by B. Riemann in 1858 who was investigating the mathematical correctness of the problem of the propagation and superposition of waves admitted by this fluid dynamics system of PDEs \cite{Riemann}. Since then, the problem of the superposition of waves has been investigated by many authors (see e.g. \cite{CourantHilbert,Jeffrey,John,Majda,PeradzynskiGeo,Rozdestvenski,Whitham}). A number of attempts to generalize the Riemann invariants method and its various applications can be found in the recent literature of the subject (see e.g. \cite{DubrovinNovikovHamil,DubrovinNovikovHydro,Tsarev,MokhovFerapontov,Ferapontov,JohnKlainerman,GrundlandWinternitz} and references therein).
	
	The results obtained from the conditional symmetry method and the generalized method of characteristics for constructing multiple Riemann wave solutions are so promising that it seems to be worthwhile to compare these methods and check their effectivness for the case of Riemann $k$-wave solutions. This is, in short, the aim of this paper.

\section{Rank 1 solutions}

Consider a properly determined homogeneous hyperbolic system of $q$ PDEs in $p$ independent variables written in the matrix form
	\begin{align}\label{rk1_sol_syst}
		\begin{aligned}
			&A^i(u)u_i=0,\quad u_i=\dfrac{\partial u}{\partial x^i},\quad i=1,...,p\\
			&x=\left(x^1,...,x^p\right)\in\mathscr{D}\subset\mathbb{R}^p,\quad u=\left(u^1,...,u^q\right)\in\mathscr{U}\subset\mathbb{R}^q,
		\end{aligned}
	\end{align}
	where $A^1,...,A^p$ are $q\times q$ matrix functions of $u$. Here, we use the summation convention. All considerations here are local. It suffices to search for solutions defined on a neighborhood of $x=0$. All functions are assumed to be smooth. The form of the system (\ref{rk1_sol_syst}) is invariant under linear transformations of the independent variables and arbitrary transformations of the dependent variables.
	
	Let us classify the solutions of the system (\ref{rk1_sol_syst}) according to their rank. For example, $u(x)$ has rank $0$ if and only if $u$ is constant and it is a trivial solution of (\ref{rk1_sol_syst}). This fact is needed for the study of rank one solutions and of more general cases of multi-wave solutions.
	
	Let us consider a wave vector which is a non-zero vector function $\lambda(u)=\left(\lambda_1(u),...,\lambda_p(u)\right)$ such that $\ker\left(\lambda_iA^i\right)\neq0$. The scalar function $r:\mathbb{R}^{p+q}\rightarrow\mathbb{R}$ given by
	\begin{equation}
		r(x,u)=\lambda_i(u)x^i	
	\end{equation}	
	is called the Riemann invariant associated with the wave vector $\lambda$. The implicit equation
	\begin{equation}\label{rk1_sol_implicit_u}
		u=f\left(r(x,u)\right)
	\end{equation}
	defines a unique function $u(x)$ on a neighborhood of $x=0$ for any function $f : \mathbb{R}\rightarrow\mathbb{R}^q$ and the matrix of derivatives of $u$ has the factorized form
	\begin{equation}
	\label{Eqq}
		\dfrac{\partial u^\alpha}{\partial x^i}(x)=\phi(x)^{-1}\lambda_i\left(u(x)\right)f'^\alpha\left(r(x,u)\right),\quad f'^\alpha=\dfrac{df^\alpha}{dr},
	\end{equation}
	where the scalar function $\phi(x)$ is given by
	\begin{equation}
		\phi(x)=\mathbbm{1}-\dfrac{\partial r}{\partial u^\alpha}(x,u(x))f'^\alpha(r(x,u(x)))\neq0,
	\end{equation}
	for which we assume that \eqref{Eqq} does not admit a gradient catastrophe. Note that $u(x)$ has rank at most equal to 1 $\left(\text{i.e. }\rank\dfrac{\partial u^\alpha}{\partial x}\leq1\right)$.
	
	If the $(p-1)$ vector fields
	\begin{equation}
		\xi_a(u)=\left(\xi_a^1(u),...,\xi_a^p(u)\right)\quad a=1,..,p-1
	\end{equation}
	satisfy the orthogonality conditions
	\begin{equation}
		\lambda_i\xi_a^i=0,
	\end{equation}
	(i.e. the vector fields $\lambda,\xi_1,\ldots,\xi_{p-1}$ form a base in $\mathbb{R}^p$), then
	\begin{equation*}
	\xi^i_a(u(x))\dfrac{\partial u^{\alpha}}{\partial x^i}(x)=0.
	\end{equation*}
	so $u^1(x),\ldots,u^q(x)$ are invariants of the vector fields $\xi^i_a(u(x))\dfrac{\partial}{\partial x^i}$ on $\mathbb{R}^p$.
	Hence, the $p$-dimensional submanifold $\{(x, u(x))\}$ is invariant under the vector field
	\begin{equation}\label{rk1_sol_vector_field}
		X_a=\xi_a^i(u)\dfrac{\partial}{\partial x^i}\quad a=1,\ldots,p-1,
	\end{equation}
	on $\mathbb{R}^p\times\mathbb{R}^q$ space. Conversely, if $u(x)$ is a $q$-component function on a neighborhood of $x=0$ such that the graph of the solution $\{(x,u(x))\}$ is invariant under all vector fields (\ref{rk1_sol_vector_field}), where $\lambda_i\xi_a^i=0$, then $u(x)$ is the solution of (\ref{rk1_sol_implicit_u}) for some function $f$. This geometrically characterizes the solutions of the implicit equation (\ref{rk1_sol_implicit_u}). Now, if $u(x)$ is the solution of (\ref{rk1_sol_implicit_u}), then we get
	\begin{equation}
		A^i(u(x))u_i=\phi(x)^{-1}\lambda_i(u(x))A^i(u(x))f'(r(x,u(x))).
	\end{equation}
	Thus $u(x)$ is a solution of the initial system (\ref{rk1_sol_syst}) if and only if the system of first-order ordinary differential equations (ODEs) for $f$
	\begin{equation}\label{rk1_sol_sol_condition}
		\lambda_i(f)A^i(f)f'=0,\quad f'=\dfrac{df}{dr},
	\end{equation}
	is satisfied, i.e. if and only if $f'$ takes values in $\ker(\lambda_iA^i)$. Note that (\ref{rk1_sol_sol_condition}) is an underdetermined system of first-order ordinary differential equations (ODEs) for $f$. The extent to which (\ref{rk1_sol_sol_condition}) constrains $f$ depends on the size of $\ker(\lambda_iA^i)$. For example, if $\lambda_iA^i=0$, then there is no constraint on $f$ at all. In any case, this exhibits a class of solutions $u(x)$ of (\ref{rk1_sol_syst}). The rank of $u(x)$ is at most equal to 1 and the graph of a solution \{(x,u(x))\} is invariant under the vector fields (\ref{rk1_sol_vector_field}), where $\lambda_i\xi_a^i=0$. Note that rescaling the wave vector $\lambda$ produces the same solutions. Note also that vector fields of the form (\ref{rk1_sol_vector_field}) commute, i.e. form an Abelian distribution.
	
	\noindent { \bf Example:}\hspace{2mm} Consider the system
	\begin{equation}
	u_t+A(u)u_x=0
	\label{2.11}
	\end{equation}
	in two independent variables $t$ and $x$, and $q$ dependent variables $u^1,\ldots,u^q$. The wave vectors are the nonzero multiples of 
	\begin{equation*}
	\lambda=(-\alpha(u),1),
	\end{equation*}
where $\alpha(u)$ is an eigenvalue function of $u$ associated with the $q\times q$ matrix function $A$. A function $u(t,x)$ defined on a neighborhood of $t=0$ satisfies an equation
\begin{equation*}
u=f(x-\alpha(u)t)
\end{equation*}	
	for some $f:\mathbb{R}\rightarrow\mathbb{R}^q$ if and only if the graph of the solution $\{(t,x,u(t,x))\}$ is invariant under the vector field
	\begin{equation}
	X=\dfrac{\partial}{\partial t}+\alpha(u)\dfrac{\partial}{\partial x},
	\label{2.12}
	\end{equation}
and such a function $u(t,x)$ is a solution of \eqref{2.11} if and only if the system of ODEs
\begin{equation}
A(f)f'=\alpha(f)f',\quad f'=\dfrac{df}{d\bar{x}},
\label{2.13}
\end{equation}
	holds. The functions
	\begin{equation*}
	\bar{t}=t,\quad\bar{x}=x-\alpha(u)t,\quad\bar{u}^1=u^1,\ldots,\bar{u}^q=u^q
	\end{equation*}
	are coordinates on $\mathbb{R}^2\times\mathbb{R}^q$. It allows us to rectify the vector fields \eqref{2.12}
	\begin{equation*}
	X=\dfrac{\partial}{\partial\bar{t}}.
	\end{equation*}
	The transverse surfaces invariant under $X$ are defined by equations of the form
\begin{equation*}
\bar{u}=f(\bar{x}),
\end{equation*}	
	where $f:\mathbb{R}\rightarrow\mathbb{R}^q$ is arbitrary, i.e. $\bar{u}(\bar{t},\bar{x})=f(\bar{x})$ is the general solution of the invariant condition
	\begin{equation*}
\bar{u}_{\bar{t}}=0.
\end{equation*}	
	Augmenting the system \eqref{2.11} with the invariant condition produces the overdetermined system
	\begin{equation*}
(A(\bar{u})-\alpha(\bar{u})\mathbb{I})\bar{u}_{\bar{x}}=0,\quad \bar{u}_{\bar{t}}=0,
\end{equation*}	
with general solution
\begin{equation*}
\bar{u}(\bar{t},\bar{x})=f(\bar{x}),
\end{equation*}	
	where $f:\mathbb{R}\rightarrow\mathbb{R}^q$ satisfies \eqref{2.13}, which reproduces the classical result.

\section{Generalized method of characteristics}

The generalized method of characteristics \cite{Burnat,PeradzynskiNonlin} for constructing multiple Riemann wave solutions of quasilinear systems (\ref{rk1_sol_syst}) can be summarized as follows :
	
	Consider an overdetermined system composed of a system (\ref{rk1_sol_syst}) which is subjected to differential constraints such that all first-order partial derivatives of the unknown functions $u^\alpha$ are decomposable (with scalar functions $\xi^s(x)\neq0$)
	\begin{equation}\label{Gen_met_syst}
		\dfrac{\partial u^\alpha}{\partial x^i}=\sum_{s=1}^{k}\xi^s(x)\gamma_s^\alpha(u)\lambda_i^s(u),
	\end{equation}
	where we assume that there exist $kp$-valued real wave functions $\lambda_i^s : \mathbb{R}^q\rightarrow\mathbb{R}$ and $kq$-valued real functions $\gamma_s^\alpha : \mathbb{R}^q\rightarrow\mathbb{R}$ which satisfy the algebraic equation (called the wave relation)
	\begin{equation}\label{Gen_met_wave_relation}
		\left(A_\alpha^{i\beta}(u)\lambda_i^s\right)\gamma_s^\alpha=0,\quad s=1,...,k,\quad\beta=1,...,q.
	\end{equation}
	Note that $u(x)$ given by (\ref{Gen_met_syst}) has rank equal to $k$ $\left(\text{i.e. }\rank\dfrac{\partial u^\alpha}{\partial x^i}\leq k\right)$. We require that every triple wave vector $\lambda^{s_1}, \lambda^{s_2}, \lambda^{s_3}$ for $s_1<s_2<s_3=1,...,k$ be linearly independent. We also assume that we have $k$ linearly independent vector functions $\gamma_s$ and we identify them with vector fields on $\mathscr{U}$-space. The algebraization of the PDEs (\ref{rk1_sol_syst}), given by (3.1) and (3.2), allows us to construct a more general class of solutions than the one given by (2.3). We have postulated a form of solution $u(x)$ for which the tangent map $du(x)$ consists of the linear combinations of simple elements $\gamma_s\otimes\lambda^s$, with non-vanishing functions $\xi^s(x)$ at any point $x\in\mathscr{D}$, namely 
	\begin{equation}\label{Gen_met_one-forms}
		du^\alpha=\sum_{s=1}^{k}\xi^s(x)\gamma_s^\alpha(u)\lambda_i^s(u)dx^i.
	\end{equation}
	We verify whether the system of one-forms (\ref{Gen_met_one-forms}) is in involution for any function $\xi^s(x)\neq0$. Here $u^1,...,u^q, \xi^1,...,\xi^k$ are considered to be unknown functions of $x^1,...,x^p$. The involutivity conditions in the Cartan sense \cite{Cartan} for the tangent map $du$ require that the commutators for each pair of vector fields $\gamma_i$ and $\gamma_j$ are spanned by these fields \cite{PeradzynskiGeo}
	\begin{equation}
		[\gamma_i,\gamma_j]\in\text{span }\{\gamma_i,\gamma_j\},\quad i\neq j=1,...,k
	\end{equation}
	and that the Lie derivatives of the wave vector $\lambda^i$ along the vector field $\gamma_j$ be
	\begin{equation}
		\mathscr{L}_{\gamma_j}\lambda^i\in\text{ span }\{\lambda^i,\lambda^j\},\quad i\neq j=1,...,k.
	\end{equation}
The conditions (3.4) and (3.5) are necessary and sufficient conditions for the existence of solutions of the one-forms (3.3) [14].	
	In practical applications, it is recommended to choose a holonomic system for the vector fields $\{\gamma_1,...,\gamma_k\}$ by requiring a proper length for each pair of vector fields $\gamma_i$ and $\gamma_j$ such that the rescaled vector fields $\tilde{\gamma_1},...,\tilde{\gamma_k}$ form an Abelian distribution.
	\begin{equation}
		[\tilde{\gamma_i},\tilde{\gamma_j}]=0,\quad i\neq j=1,...,k.
	\end{equation}
	Here, we are interested in rescaling the vector fields $\{\gamma_1,...,\gamma_k\}$ by functions of $u$, which ensures that the rescaled vector fields still satisfy (\ref{Gen_met_wave_relation}). Note that the Frobenius theorem is not sufficient to prove the existence of an appropriate rescaling of these vector fields $\{\gamma_1,...,\gamma_k\}$ leading them to commute among themselves. The existence of $k$-wave solutions requires a more restrictive condition on the form of the commutator between the vector fields $\gamma_i$ and $\gamma_j$ than the one assumed by the Frobenius theorem (see e.g. \cite{Stenberg}). This fact requires a modification of this theorem for integration which is useful for the purpose of constructing and investigating $k$-wave superpositions expressed in terms of Riemann invariants in $\mathbb{R}^p$-space.
	
	\textbf{Theorem 1.} (the modified Frobenius theorem by rescaling \cite{GrundlandDeLucas}.)
	
	Suppose that a family of vector fields $\gamma_1,...,\gamma_k$ is defined on a $q$-dimensional manifold $N$ such that
	\begin{align}
		\begin{aligned}
			i)&\ \gamma_1\wedge...\wedge\gamma_k\neq0,\\
			ii)&\ [\gamma_i,\gamma_j]=h_{ij}^i\gamma_i+h_{ij}^j\gamma_j, \text{ (no summation convention)}
		\end{aligned}
	\end{align}
	for certain functions $h_{ij}^i\in\mathbb{C}^\infty(N)$ with $i,j,l=1,...,k$. Then there exists a set of functions $f_1,...,f_k\in\mathbb{C}^\infty(N)$ such that the rescaled vector fields $\gamma_1,...,\gamma_k$ form an Abelian distribution
	\begin{equation}\label{Gen_met_Abelian_distribution}
		[f_i\gamma_i,f_j\gamma_j]=0,\quad i,j=1,...,k
	\end{equation}
	for all $i,j=1,...,k$. 
	
	The proof uses the induction hypothesis and  the above result is demonstrated in detail in \cite{GrundlandDeLucas}.

	Suppose now, according to Theorem 1, that the normalized $k$-dimensional distribution $\{\gamma_1,...,\gamma_k\}$ is Abelian (i.e. that the conditions (\ref{Gen_met_Abelian_distribution}) hold). This means that we can determine a $k$-dimensional manifold $S$ on $\mathbb{R}^q$ obtained by integrating the nonlinear system of $kq$ PDEs
	\begin{equation}\label{Gen_met_integration_kq_PDEs}
		\dfrac{\partial f^\alpha}{\partial r^s}=\gamma_s^\alpha(f^1,...,f^q),\quad s=1,...,k
	\end{equation}
	with a solution defined by 
	\begin{equation}
	\label{S}
		S:u=\left(f^1(r^1,...,r^k),...,f^q(r^1,...,r^k)\right).
	\end{equation}
	Hence the vector fields $\gamma_s=\gamma_s^\alpha\dfrac{\partial}{\partial u^\alpha}$ are Frobenius integrable and, due to the independence condition (3.7$i$), these vector fields determine a foliation of $\mathbb{R}^q$ given by $k$-dimensional leaves.	Note that a solution of the system of one-forms (3.3) also satisfies the original system of PDEs (\ref{rk1_sol_syst}). Hence the system (3.3) can be parametrized in terms of the variables $r^1,...,r^k$ if $S$ admits a coordinate system \eqref{S}. The matrix of derivatives of (3.10) takes the form
	\begin{equation}
		\dfrac{\partial u^\alpha}{\partial x^i}=\sum_{s=1}^{k}\dfrac{\partial u^\alpha}{\partial r^s}\dfrac{\partial r^s}{\partial x^i}.
	\end{equation}
	Assume that the wave functions $\lambda_i^s(u)$ are pulled back to the manifold $S\subset\mathscr{U}$. Then the $kp$ wave functions $\lambda_i^s(u)$ become functions of a coordinate system $r^1,...,r^k$ on $S$.	From the linear independence of the vector fields $\gamma_1,...,\gamma_k$, and comparing (3.3) with (3.11), we obtain the Pfaffian system of one-forms 
	\begin{equation}
		dr^s=\xi^s(x)\lambda_i^s(r^1,...,r^k)dx^i,\quad s=1,...,k\quad(\text{with } \xi^s(x)\neq0).
	\end{equation}
	Note that the elimination of the functions $\xi^1,...,\xi^k$ in the Pfaffian system (3.12) leads to an overdetermined system of PDEs. Hence, it is convenient from the computational point of view to impose on the wave vectors $\lambda^1,...\lambda^k$ the involutivity conditions which guarantee the existence of a certain class of solutions of the Pfaffian system (3.12) with the freedom of $k$ arbitrary functions of one variable. Under the above assumption, the conditions (3.5) on the one-forms $\lambda^s$ can be written in terms of $r^1,...,r^k$ and has to satisfy the linear first-order system of PDEs \cite{PeradzynskiNonlin}
	\begin{equation}\label{Gen_met_1st_order_PDEs}
		\dfrac{\partial \lambda^s}{\partial r^p}=\alpha_p^s\lambda^s+\beta_p^s\lambda^p\quad\text{ (no summation)}
	\end{equation}
	for all $s\neq p=1,...,k$, where $\alpha_p^s$ and $\beta_p^s$ are functions of $r^1,...,r^k$. Next, after an appropriate normalization of each wave vector $\lambda^s$ such that $\lambda^s=\left(1,\lambda_2^s,...,\lambda_p^s\right)$, the involutivity conditions (3.13) can be written
	\begin{equation}
		\dfrac{\partial\lambda_\mu^s}{\partial r^l}=\alpha_l^s\left(\lambda_\mu^s-\lambda_\mu^l\right),\quad \dfrac{\partial\lambda_\mu^l}{\partial r^s}=\beta_s^l\left(\lambda_\mu^l-\lambda_\mu^s\right),\quad s\neq l,\quad \mu=1,...,p.
	\end{equation}
	Cross-differentiating the equations (3.14), we obtain a linear hyperbolic system of $kp$ second-order equations for each component of the wave function $\lambda_\mu^s$, that is \cite{GrundlandVassiliou}
	\begin{equation}
		\dfrac{\partial^2\lambda_\mu^s}{\partial r^s\partial r^l}+\tilde{\alpha}_l^s\dfrac{\partial\lambda_\mu^s}{\partial r^{(s)}}+\tilde{\beta}_l^s\dfrac{\partial\lambda_\mu^s}{\partial r^{(l)}}=0,\quad \mu=1,...,p
	\end{equation}
	for $s\neq l=1,...,k$ (no summation over the indices $s$ and $l$), where $\tilde{\alpha}_l^s$ and $\tilde{\beta}_l^s$ are functions of $r^1,..,r^k$. Equation (3.15) always admits a solution by the method of characteristics (it is a Darboux problem). We look for the most general solution of the linear second-order PDEs (3.15) which allows us in turn to determine the solution of the Pfaffian system (3.12) which has the implicit form \cite{PeradzynskiGeo,GrundlandVassiliou}
	\begin{equation}
		\lambda_i^s(r^1,...,r^k)x^i=\psi^s(r^1,...,r^k).
	\end{equation}
	The quantities $r^1,...,r^k$ constitute Riemann invariants as they are implicitly defined as functions of $x^1,...,x^p$ by the relation (3.16) in the space $\mathbb{R}^p$. Thus, the construction of the $k$-wave solution is determined by the implicit relations
	\begin{align}
		\begin{aligned}
			&u=\left(f^1(r^1,...,r^k),...,f^q(r^1,...,r^k)\right),\\
			&\lambda_i^s(r^1,...,r^k)x^i=\psi^s(r^1,...,r^k),\quad s=1,...,k.
		\end{aligned}
	\end{align}
	This completes the factorization of the problem of constructing $k$-wave solutions by the generalized method of characteristics. This approach has produced many new analytic solutions of hydrodynamic-type systems (equations appearing in continuous media, e.g. classical and relativistic hydrodynamics (see e.g.\cite{Zajaczkowski,Lichnerowicz,GrundlandZelazny}  and references therein) and in nonlinear field equations (see e.g. \cite{Boillat,Rendal,Trautman} and references therein)).

 \section{The elastic superposition of two waves}
 
 Let us now demonstrate that certain classes of solutions of the Pfaffian system (3.12) describe the elastic superposition of two waves. For the sake of simplicity this statement is illustrated for a quasilinear system with two independent variables (time $t$ and one space variable $x$)
	\begin{equation}
		u_t+A(u)u_x=0,
	\end{equation}
	where $A$ is a $q\times q$ matrix function of $u=\left(u^1,...,u^q\right)$. We assume that the double wave  solutions of (3.3) exist (for $k=2$) and that, in the initial conditions for time $t=t_0$, the non-vanishing dependent variables, $\xi^1(t_0,x)$ and $\xi^2(t_0,x)$, are different from zero.
	
	It was proved in \cite{CourantFriedrichs,John} that if the initial data are sufficiently small at $t=t_0$, then there exists a time interval $[t_0,T]$ in which the gradient catastrophe for the derivatives of the solutions $\dfrac{\partial r^s}{\partial x},\ s=1,2$, of the Pfaffian system (3.12) does not occur. Each function $r^s(t,x),\ s=1,2$ is constant along the appropriate family of characteristics
	\begin{equation}
		C^{(s)} : \dfrac{dx}{dt}=\nu_s\left(r^1(t,x),r^2(t,x)\right),\quad s=1,2
	\end{equation} 
	of the system obtained from the Pfaffian system (3.12) by the elimination of the variables $\xi^s$
	\begin{equation}
		\dfrac{\partial r^s}{\partial t}+\nu_s\left(r^1(t,x),r^2(t,x)\right)\dfrac{\partial r^s}{\partial x}=0,
	\end{equation}
	where the wave vectors $\lambda^s$ have been normalized to $\lambda^s=\left(\nu_s(r^1,r^2),-1\right)$ and $\nu_s$ are the eigenvalues associated with the matrix $A$. If we choose, in the space of independent variables $\mathbb{R}^2$, the initial data at $t=t_0$ in such a way that the derivatives $\dfrac{\partial r^s}{\partial x}$ have compact and disjoint supports
	\begin{align}
		\begin{aligned}
			&1)\quad t=t_0 : \text{supp }\dfrac{\partial r^s}{\partial x}(t_0,x)\subset[a_s,b_s],\\
			&2)\quad t=t_0 : \text{supp }\dfrac{\partial r^1}{\partial x}(t_0,x)\cap\text{supp }\dfrac{\partial r^2}{\partial x}(t_0,x)=\phi,
		\end{aligned}
	\end{align}
	where $a_s,b_s\in\mathbb{R},\quad s=1,2$, and $a_1<b_1<a_2<b_2$, then, for an arbitrary time $t_0<t<T$, the supports $\text{supp }\dfrac{\partial r^s}{\partial x}(t,x)$ are contained in the strips between appropriate characteristics of the families $C^{(s)}$ passing through the ends of the intervals $[a_s,b_s]$. It was proved \cite{John,JohnKlainerman} that the initial data at $t=t_0$ are sufficiently small and also satisfy the eigenvalue condition
	\begin{equation}
		\exists\ c>0\ \forall\ (t_0,x),(t_0,x')\in[t_0,T]\times\mathbb{R} : \nu_1(r(t_0,x))-\nu_2(r(t_0,x'))\geq c,
	\end{equation}
	for any $x'<x$, where $r=\left(r^1,r^2\right)$. This means that every characteristic of the same family ($C^{(1)}$) has an inclination (measured with respect to the positive direction of the $x$ axis) which is smaller than the slope of any charateristic of the other family $C^{(2)}$. In this case the strips containing $\text{supp }\dfrac{\partial r^s}{\partial x}(t,x)$ divide the remaining part of the space of independent variables $\mathbb{R}^2$ into four disjoint regions, say I-IV respectively. 
	In the region $B : \mathbb{R}^2\smallsetminus\{\text{supp }\dfrac{\partial r^1}{\partial x}(t,\cdot)\cup\text{supp }\dfrac{\partial r^2}{\partial x}(t,\cdot)\}$	the solution $r^s(t,x)$ of the Pfaffian system (3.12) is constant i.e. $r^s=r_0^s.$
		
	Let us denote by $t_1$ and $t_2$ the times in which $\text{supp }\dfrac{\partial r^1}{\partial x}(t,x)$ and $\text{supp }\dfrac{\partial r^2}{\partial x}(t,x)$ have only one common point. For times $t\in [t_0,t_1]$ we have
	\begin{equation}
		\text{supp }\dfrac{\partial r^1}{\partial x}(t,x)\cap\text{supp }\dfrac{\partial r^2}{\partial x}(t,x)=\phi,\quad t<t_1
	\end{equation}
	Therefore, the solution can be interpreted as the propagation of two separate (non-interacting) simple waves. For the times $t\in[t_1,t_2]$, the characteristics of the families $C^{(1)}$ and $C^{(2)}$ containing $\text{supp }\dfrac{\partial r^s}{\partial x}(t,x)\ s=1,2$ cross each other, i.e.
	\begin{equation}
		\text{supp }\dfrac{\partial r^1}{\partial x}(t,x)\cap\text{supp }\dfrac{\partial r^2}{\partial x}(t,x)\neq\phi,\quad t\in[t_1,t_2]
	\end{equation}
	We interpret this phenomenon as a superposition of two simple waves. They constitute a double wave solution of (3.12). For time $t\in[t_2,T]$, making use of the conditions (4.4) and (4.5), the strips containing the supports of the simple waves separate again. So, we again have the situation for which (4.6) holds. This means that the double wave solution decays in an exact way into two separate simple waves of the same type as imposed in the initial data for $t=t_0$. Therefore, the number and type of simple waves is conserved. Therefore, we can speak of the elastic superposition of simple waves.
	
	In general, in the case of more than two simple waves $(k>2)$, the superposition of waves is analogous but more involved, because the region $B\subset\mathbb{R}^2$ is divided by the supports $\text{supp }\dfrac{\partial r^s}{\partial x},\ s=1,...,k$ into $2^k$ subregions.

\section{Conditional symmetry method} 

There now arises the question of what additional results, if any, can be achieved by applying symmetry group analysis to the problem of the construction of Riemann $k$-waves. Let us consider a fixed set of $k$ linearly independent wave vectors $\lambda^1,...,\lambda^k,\ 1\leq k\leq p$ with corresponding Riemann invariants
	\begin{equation}
		r^1(x,u)=\lambda_i^1(u)x^i,...,r^k(x,u)=\lambda_i^k(u)x^i.
	\end{equation}
	The implicit equation
	\begin{equation}
		u=f\left(r^1(x,u),...,r^k(x,u)\right)
	\end{equation}
	defines a unique function $u(x)$ on a neighborhood of $x=0$, for any function $f : \mathbb{R}^k\rightarrow\mathbb{R}^q$. For all $i$ and $\alpha$, the decomposition of the matrix of derivatives of $u$ has the factorized form
	\begin{equation}
		\dfrac{\partial u^\alpha}{\partial x^i}(x)=\left(\phi(x)^{-1}\right)_j^l\lambda_i^j(u(x))\dfrac{\partial f^\alpha}{\partial r^l}(r(x,u(x))),
	\end{equation}
	where $\phi(x)$ is an invertible $k\times k$ matrix function of $x$ (we assumed no gradient catastrophe)
	\begin{equation}
		\phi(x)_j^i=\delta_j^i-\dfrac{\partial r^i}{\partial u^\alpha}(x,u(x))\dfrac{\partial f^\alpha}{\partial r^j}(f(x,u(x))),\quad\dfrac{\partial r^i}{\partial u^\alpha}=\dfrac{\partial\lambda_j^i}{\partial u^\alpha}x^j.
	\end{equation}
	Note that the form of the matrix of derivatives (5.3) is more general than the one obtained from the method of characteristics (3.1), where we assumed that $\xi^i(x)$ are scalar functions. The rank of $u(x)$, given by (5.3), is at most equal to $k$ $\left(\text{i.e. }\rank\dfrac{\partial u^\alpha}{\partial x^i}\leq k\right)$. If the $(p-k)$ vector fields
	\begin{equation}
		\xi_a(u)=\left(\xi_a^1(u),...,\xi_a^p(u)\right)^T,\quad a=1,...,p-k
	\end{equation}
	satisfy the orthogonality conditions
	\begin{equation}\label{Cond_symm_orth_cond}
		\lambda_i^j\xi_a^i=0,\quad j=1,...,k
	\end{equation}
	then the $p$-dimensional manifold $\{(x,u(x))\}$ is invariant under the $(p-k)$ vector fields
	\begin{equation}\label{Cond_symm_vector_field}
		X_a=\xi_a^i(u)\dfrac{\partial}{\partial x^i},\quad a=1,...,p-k
	\end{equation}
	on $\mathbb{R}^p\times\mathbb{R}^q$. Conversely, if $u(x)$ is a $q$-component function on a neighborhood of $x=0$ such that the graph $\{(x,u(x))\}$ is invariant under the $(p-k)$ vector fields (5.7) for which (5.6) holds, then $u(x)$ is the solution of (5.2) for some $f$ because $\{r^q,...,r^k,u^1,...,u^q\}$ constitutes a complete set of invariants of the Abelian algebra of vector fields (5.7). If $u(x)$ is the solution of the implicit equation (5.2), then the system (\ref{rk1_sol_syst}) reduces to the system of q PDEs written in the matrix form
	\begin{equation}
		A^i(u(x))\dfrac{\partial u^\alpha}{\partial x^i}(x)=\left(\phi(x)^{-1}\right)_j^l\lambda_i^j(u(x))A^i(u(x))\dfrac{\partial f^\alpha}{\partial r^l}(r(x,u(x)))=0,
	\end{equation}
	where the matrix function $\phi(x)$ is given by (5.4). Comparing the result (5.3) with (3.1), we observe that we cannot generally expect a reduction analogous to (5.8) because $\phi$ is no longer a scalar function.
	
	A change of variables illuminates the construction of $k$-wave solutions of (\ref{rk1_sol_syst}). Let us consider a coordinate transformation and let us fix $k$ linearly independent wave vectors $\lambda^1,...,\lambda^k$. Consider the $p\times k$ matrix whose columns are the vectors $\lambda^1,...,\lambda^k$ and then assume that there exists a $k\times k$ invertible submatrix function of $u$ given by
	\begin{equation}
		\Lambda(u)=\left(\lambda_l^j(u)\right).\quad i,j=1,...,k
	\end{equation}
	Then the $(p-k)$ independent vector fields
	\begin{equation}
		X_a=\dfrac{\partial}{\partial x^a}-\sum_{l,j=1}^{k}\left(\Lambda(u)^{-1}\right)_j^l\lambda_a^j\dfrac{\partial}{\partial x^l},\quad a=k+1,...,p
	\end{equation}
	have the form (5.7) with the orthogonality property (5.6). The sets of functions
	\begin{equation}\label{Cond_symm_var_change}
		\bar{x}^1=r^1(x,u),...,\bar{x}^k=r^k(x,u),\bar{x}^{k+1}=x^{k+1},...,\bar{x}^p=x^p,\bar{u}=u^1,...,\bar{u}^q=u^q,
	\end{equation}
	are coordinates on $\mathbb{R}^p\times\mathbb{R}^q$. This allows us to rectify the vector fields (5.10)
	\begin{equation}
		X_{k+1}=\dfrac{\partial}{\partial\bar{x}^{k+1}},...,X_{p}=\dfrac{\partial}{\partial\bar{x}^{p}}.
	\end{equation}
	The $p$-dimensional manifold $\{(x,u(x))\}$ is invariant under (5.12) and is defined by the equation
	\begin{equation}
		\bar{u}=f\left(\bar{x}^1,...,\bar{x}^k\right),
	\end{equation}
	where $f : \mathbb{R}^k\rightarrow\mathbb{R}^q$ is an arbitrary function of $\bar{x}^1,...,\bar{x}^k$. This means that the expression (5.13) is the general solution of the invariance conditions
	\begin{equation}
		\bar{u}_{k+1}=...=\bar{u}_p=0.
	\end{equation}
	The initial system of PDEs (\ref{rk1_sol_syst}) can be described in terms of the new coordinates (5.11), i.e. in terms of $\bar{x}$, $\bar{u}$ by
	\begin{equation}
		\bar{A}^i(\bar{x},\bar{u},\bar{u}_{\bar{x}})\bar{u}_i=0,
	\end{equation}
	where
	\begin{equation*}
		\bar{A}^1=\dfrac{Dr^1}{Dx^i}A^i,...,\bar{A}^k=\dfrac{Dr^k}{Dx^i}A^i,\bar{A}^{k+1},...,\bar{A}^p=A^p
	\end{equation*}
	and
	\begin{equation*}
		\dfrac{Dr^i}{Dx^j}|_{\bar{u}^{k+1}=...=\bar{u}^p=0}=\left(\phi^{-1}\right)_l^i\lambda_j^l,\quad\phi_j^i=\delta_j^i-\dfrac{\partial r^i}{\partial u^\alpha}\bar{u}_j^\alpha.
	\end{equation*}
	So, subjecting (5.15) to the invariant conditions (5.14), we obtain the overdetermined system of PDEs which defines Riemann $k$-wave solutions of the initial system (\ref{rk1_sol_syst})
	\begin{align}
		\begin{aligned}
			\sum_{i,j=1}^{k}\sum_{l=1}^{p}\left(\phi^{-1}\right)_j^i\lambda_l^jA^l\bar{u}_i=0,\quad\text{where }\left(\phi^{-1}\right)_j^i=\delta_j^i-\dfrac{\partial r^i}{\partial u^\alpha}\bar{u}_j^\alpha,\\ \bar{u}_{k+1}=...=\bar{u}_p=0,\quad i,j=1,...,k.
		\end{aligned}
	\end{align}
	As one can observe, this approach is a variant of the conditional symmetry method.
	
	Numerous examples of new Riemann $k$-wave solutions of hydrodynamic-type systems in $(n+1)$ dimensions have been obtained through this approach. This approach was found to have an advantage over the method of characteristics, resulting mainly from the relaxing of the second-order differential constraints (3.15) for the wave functions $\lambda^s$. It delivers larger classes of solutions expressed in terms of Riemann invariants through a simpler procedure (see e.g. \cite{GrundlandTafel,GrundlandHuard,DoyleGrundland,GrundlandLamothe} and references therein).

\section{Examples} 

We provide some examples to illustrate the construction of solutions which are expressible in terms of Riemann invariants.

\subsection{Barotropic fluid flow in $(n+1)$ dimensions}

We present a simple example to illustrate the construction introduced in section 5. Consider the barotropic system (i.e. a properly determined hyperbolic system of PDEs)
	\begin{equation}
		u_t+(u\cdot\nabla)u=0,\quad\rho_t+\nabla(\rho u)=0,\quad\rho>0,
	\end{equation}
	which describes an inviscid fluid flow at constant pressure in the $(n+1)$-dimensional space. These equations represent momentum and mass conservation. There are $(n+1)$ dependent variables $\left(u^1,...,u^n,\rho\right)$. The functions
	\begin{equation*}
		\lambda^1=\left(-u^1,1,0,...,0\right),...,\lambda^n=\left(-u^n,0,...,0,1\right)
	\end{equation*} 
	are wave vectors and the corresponding vector fields
	\begin{equation}
		X_1=\dfrac{\partial}{\partial t}+u^1\dfrac{\partial}{\partial x^1},...,X_n=\dfrac{\partial}{\partial t}+u^n\dfrac{\partial}{\partial x^n}
	\end{equation}
	span an Abelian distribution $\mathscr{A}$. Augmenting (6.1) with $\mathscr{A}$-invariance conditions produces the overdetermined system
	\begin{equation}
		\nabla\cdot u=0,\quad u_t+(u\cdot\nabla)u=0,\quad\rho_t+(u\cdot\nabla)\rho=0,\quad\rho>0.
	\end{equation}
	Thus, a solution of (6.1) is $\mathscr{A}$-invariant if and only if the velocity field $\mathscr{u}$ is divergence-free. We will now derive the general solution of (6.1), as well as the class of $\mathscr{A}$-invariant solutions. The system (6.1) is described in the coordinates
	\begin{equation}
		\bar{t}=t,\quad\bar{x}^1=x^1-u^1t,\ldots,\bar{x}^n=x^n-u^nt,\quad\bar{u}^1=u^1,\ldots,\bar{u}^n=u^n,\quad\bar{\rho}=\rho
	\end{equation}
	by
	\begin{equation}
		\bar{u}_{\bar{t}}=0,\quad\bar{\rho}_{\bar{t}}+\bar{\rho}\tr\left(\left(\mathbbm{1}+t\bar{u}_{\bar{x}}\right)^{-1}\bar{u}_{\bar{x}}\right)=0,\quad\rho>0,\quad\text{where }\bar{u}_{\bar{x}}=\dfrac{\partial(u^1,...,u^n)}{\partial(x^1,...,x^n)}.
	\end{equation}
	The general solution of	$\bar{u}_{\bar{t}}=0$ is
	\begin{equation*}
		\bar{u}=f(\bar{x}),
	\end{equation*}
	where $f : \mathbb{R}^n\rightarrow\mathbb{R}^n$ is arbitrary. Then, using the Abel identity, we get
	\begin{equation*}
		\tr\left(\left(\mathbbm{1}+\bar{t}\dfrac{\partial\bar{u}}{\partial\bar{x}}(\bar{x})\right)^{-1}\dfrac{\partial\bar{u}}{\partial\bar{x}}(\bar{x})\right)=\dfrac{\partial}{\partial\bar{t}}\ln\left(\det\left(\mathbbm{1}+tDf(\bar{x})\right)\right),\quad\text{ where } Df(\bar{x})=\dfrac{\partial(f^1,...,f^n)}{\partial(\bar{x}^1,...,\bar{x}^n)}.
	\end{equation*}
	So the condition on $\bar{\rho}(\bar{t},\bar{x})$ is
	\begin{equation*}
		\dfrac{\partial}{\partial\bar{t}}\ln\left(\bar{\rho}\left(\bar{t},\bar{x}\right)\det\left(\mathbbm{1}+\bar{t}Df(\bar{x})\right)\right)=0.
	\end{equation*}
	Thus, the general solution $\left(u(t,x),\rho(t,x)\right)$ of (6.1) is given implicitly by \cite{DoyleGrundland}
	\begin{equation*}
		u(t,x)=f(x-ut),\quad\rho(t,x)=\dfrac{g(x-u(t,x)t)}{\det\left(\mathbbm{1}+Df(x-u(t,x)t)\right)},
	\end{equation*}
	where $g:\mathbb{R}^n\rightarrow\mathbb{R}$ is an arbitrary positive function.
	
	The $(p-n)$ vector fields (6.2) can be rectified under the coordinate transformation (6.4). Then, we get one vector field
	\begin{equation}
		X=\dfrac{\partial}{\partial\bar{t}}.
	\end{equation}
	The general solution of the invariance conditions $\bar{u}_{\bar{t}}=\bar{\rho}_{\bar{t}}=0$, is
	\begin{equation}
		\bar{u}(\bar{t},\bar{x})=f(\bar{x}),\quad \bar{\rho}(\bar{t},\bar{x})=g(\bar{x}),
	\end{equation}
	where $f : \mathbb{R}^n\rightarrow\mathbb{R}^n$, and $g : \mathbb{R}^n\rightarrow\mathbb{R}$ are arbitrary functions. Augmenting the system (6.5) with the invariance conditions (6.7) produces the system
	\begin{equation*}
		\tr\left(\left(\mathbbm{1}+\bar{t}\bar{u}_{\bar{x}}\right)^{-1}\bar{u}_{\bar{x}}\right)=0,\quad\bar{u}_{\bar{t}}=\bar{\rho}_{\bar{t}}=0,\quad\bar{\rho}>0,
	\end{equation*}
	and using the Abel identity, we get
	\begin{equation*}
		\dfrac{\partial}{\partial\bar{t}}\det\left(\mathbbm{1}+\bar{t}Df(\bar{x})\right)=0,\quad g>0.
	\end{equation*}
	This implies that the $\mathscr{A}$-invariant solutions of the barotropic flow (6.1) are described implicitly by
	\begin{equation}
		u(x,t)=f(x-ut),\quad\rho(x,t)=g(x-u(t,x)t)>0,
	\end{equation}
	where $g : \mathbb{R}^n\rightarrow\mathbb{R}$ is an arbitrary positive function and $f : \mathbb{R}^n\rightarrow\mathbb{R}^n$ has a nilpotent differential $Df$ $\left(\text{i.e. }(Df)^k=0,\quad k=1,...,n\right)$. In this case, for any degree of nilpotency of the differential $Df(x)$, the obtained equations correspond to the hypersurfaces immersed in $\mathbb{R}^{n+1}$ which have zero Gaussian curvature.

\subsection{The magnetohydrodynamic equations}

We consider an ideal compressible conductive and isothermal fluid flow in $(3+1)$-dimensions placed in the presence of a magnetic field. Under these assumptions, the system of equations forms an evolutionary hyperbolic first-order system of eight PDEs in four independent variables (time $t$ and three space variables $x=\left(x^1,x^2,x^3\right)$). The magnetohydrodynamics (m.h.d.) system takes the form \cite{Pai,Davidson}
	\begin{align}\label{mhd_system}
		\begin{split}
			\dfrac{d}{dt}\rho+\rho\nabla\cdot\bar{v}=0,\\
			\rho\dfrac{d\bar{v}}{dt}+\nabla p+\dfrac{1}{4\pi}\bar{H}\cross\left(\nabla\cross\bar{H}\right)=0,\\
			\dfrac{\partial}{\partial t}\bar{H}-\nabla\cross\left(\bar{v}\cross\bar{H}\right)=0,\quad\nabla\cdot\bar{H}=0,\\
			\dfrac{d}{dt}\left(\dfrac{p}{\rho\gamma}\right)=0,\quad\dfrac{d}{dt}=\dfrac{\partial}{\partial t}+\left(\bar{v}\cdot\nabla\right).
		\end{split}
	\end{align}
	Here, we have used the following notation : $\rho$ is the density of the fluid, $p$ is the pressure of the fluid, $\bar{v}$ is the vector field of the fluid velocity, $\bar{H}$ is the vector of the magnetic field and $\gamma$ is the polytropic exponent. We also assume that the inductive capacity $\mu$ is equal to one. The system (\ref{mhd_system}) represents the mass and momentum conservation of the field equation, the Gauss law for the magnetic field, and conservation of the entropy, respectively. The unknown functions are $u=\left(\rho,p,\bar{v},\bar{H}\right)\in\mathbb{R}^8$. The Lie subgroup analysis of the symmetry groups of the m.h.d. equations \eqref{mhd_system} was established in \cite{GL95}.
	
	In this paper, we look for Alfvénian wave solutions admitted by the m.h.d. equations (\ref{mhd_system}) involving Riemann invariants. The corresponding simple Alfvén elements $A_2(\varepsilon)$ associated with the Alfvén wave velocity $\delta|\bar{\lambda}|$  (i.e. the eigenvalue corresponding to the matrix associated with the m.h.d. system (6.9))
\begin{equation*}	
\delta|\bar{\lambda}|=\dfrac{\varepsilon\bar{H}\cdot\bar{\lambda}}{\sqrt{4\pi\rho}},\qquad \varepsilon=\pm 1,
\end{equation*}
(where $\lambda^i=(\delta|\bar{\lambda}|-\bar{v}\cdot\bar{\lambda}^i,\bar{\lambda}^i)$, $i=1,2$ are linearly independent wave vectors) take the form \cite{Zajaczkowski}
	\begin{equation}
	\label{AlfvenElements}
		\gamma_i=\left(0,0,\dfrac{\varepsilon\bar{h}_i}{\sqrt{4\pi\rho}},\bar{h}_i\right),\quad\lambda^i=\left(0,\bar{\alpha}^i\cross\bar{h}_i\right),\quad\bar{H}_i\bar{h}_i=0,\quad i=1,2,
	\end{equation}
	where $\bar{h}_i=\partial\bar{H}/\partial r^i,\quad i=1,2$ are arbitrary vectors orthogonal to $\bar{H}$ and $\bar{\alpha}^i$ are arbitrary linearly independent vectors in $\mathbb{R}^3$. These simple elements $A_2(\varepsilon)$ satisfy the following system of PDEs
\begin{equation*}
\begin{split}
\rho=\rho_0,\quad p=p_0,\quad |\bar{H}\left(\bar{x}^1,\bar{x}^2\right)|^2=\mbox{const.},\quad \rho_0,p_0\in\mathbb{R},\\ \dfrac{\partial\bar{v}}{\partial \bar{t}}=\dfrac{1}{4\pi\rho_0}(\bar{H}\cdot\nabla)\bar{H},\quad \dfrac{\partial\bar{v}}{\partial \bar{t}}=\varepsilon\dfrac{(\bar{H}\cdot\nabla)\bar{v}}{\sqrt{4\pi\rho_0}},\quad\nabla\cdot\bar{v}=0,\\ \dfrac{\partial\bar{H}}{\partial \bar{t}}=(\bar{H}\cdot\nabla)\bar{v},\quad \dfrac{\partial\bar{H}}{\partial \bar{t}}=\varepsilon\dfrac{(\bar{H}\cdot\nabla)\bar{H}}{\sqrt{4\pi\rho_0}},\quad
\nabla\cdot\bar{H}=0.
\end{split}
\end{equation*}
	
	 Using the generalized method of characteristics, we search for a special class of double wave solutions (\ref{Gen_met_one-forms}) for which the tangent map $du(x)$ is the sum of two simple Alfvén elements \eqref{AlfvenElements}
	\begin{align}
		\begin{split}
			du(x)=&\xi^1(x)\gamma_1\otimes\lambda^1+\xi^2(x)\gamma_2\otimes\lambda^2\\
			=&\xi^1(x)\left(0,0,\dfrac{\varepsilon\bar{h}_1}{\sqrt{4\pi\rho}},\bar{h}_1\right)\otimes\left(0,\bar{\alpha}^1\cross\bar{h}_1\right)\\+&\xi^2(x)\left(0,0,\dfrac{\varepsilon\bar{h}_2}{\sqrt{4\pi\rho}},\bar{h}_2\right)\otimes\left(0,\bar{\alpha}^2\cross\bar{h}_2\right),\quad\bar{H}\cdot\bar{h}_i=0,\quad i=1,2
		\end{split}
	\end{align}
	The compatibility conditions (\ref{Gen_met_Abelian_distribution}) require that the vector fields $\gamma_1$ and $\gamma_2$ form an Abelian distribution (\ref{Gen_met_Abelian_distribution}). This means that there exists a two-dimensional manifold $S$ in $\mathbb{R}^8$ obtained by integrating the system of 16 PDEs (\ref{Gen_met_integration_kq_PDEs}) taking the form
	\begin{equation}\label{mhd_16_PDEs}
		\begin{aligned}
			\dfrac{\partial\rho}{\partial r^1}&=0,\\
			\dfrac{\partial p}{\partial r^1}&=0,\\
			\dfrac{\partial\bar{v}}{\partial r^1}&=\dfrac{\varepsilon\bar{h}_1}{\sqrt{4\pi\rho}},\\
			\dfrac{\partial\bar{H}}{\partial r^1}&=\bar{h}_1,
		\end{aligned}
		\qquad\qquad
		\begin{aligned}
			\dfrac{\partial\rho}{\partial r^2}&=0,\\
			\dfrac{\partial p}{\partial r^2}&=0,\\
			\dfrac{\partial\bar{v}}{\partial r^2}&=\dfrac{\varepsilon\bar{h}_2}{\sqrt{4\pi\rho}},\\
			\dfrac{\partial\bar{H}}{\partial r^2}&=\bar{h}_2,
		\end{aligned}
		\qquad\qquad
		\bar{H}\cdot\bar{h}_i=0,\quad i=1,2
	\end{equation}
	where $r^1$ and $r^2$ are a coordinate system on $S$. The wave vector fields $\lambda_1$ and $\lambda_2$ have to satisfy the linear system (\ref{Gen_met_1st_order_PDEs}) of 8 PDEs,
	\begin{equation}\label{mhd_8_PDEs}
		\begin{aligned}
			\dfrac{\partial\lambda_0^1}{\partial r^1}&=0,\\
			\dfrac{\partial\bar{\lambda}^1}{\partial r^1}&=\beta^1\left(\bar{\alpha}^1\cross\bar{h}_1\right)+\beta^2\left(\bar{\alpha}^2\cross\bar{h}_2\right),
		\end{aligned}
		\qquad\qquad
		\begin{aligned}
			\dfrac{\partial\lambda_0^2}{\partial r^2}&=0,\\
			\dfrac{\partial\bar{\lambda}^2}{\partial r^2}&=\beta^3\left(\bar{\alpha}^1\cross\bar{h}_1\right)+\beta^4\left(\bar{\alpha}^2\cross\bar{h}_2\right),
		\end{aligned}
	\end{equation}
	where $\beta^1,...,\beta^4$ are some functions of $r^1$ and $r^2$. The double Alfvén wave solution of the system (\ref{mhd_16_PDEs}) and (\ref{mhd_8_PDEs}) written in terms of Riemann invariants takes the form \cite{Zajaczkowski}
	\begin{equation}\label{mhd_Alfven_form}
		\rho=\rho_0\quad p=p_0,\quad |\bar{H}\left(r^1,r^2\right)|^2=\mbox{const.},\quad\bar{v}=\dfrac{\varepsilon\bar{H}}{\sqrt{4\pi\rho_0}},\quad\varepsilon=\pm1,
	\end{equation}
	where $\bar{H}$ is an arbitrary vector with constant length and $\rho_0,p_0,H_0$ are some constants. The Riemann invariants $r^1$ and $r^2$ are determined from (\ref{Gen_met_1st_order_PDEs}), which leads to the following conditions on the Riemann invariants $r^1$ and $r^2$, i.e.
	\begin{equation}\label{mhd_characteristics_method}
		\dfrac{\partial\bar{H}}{\partial r^1}\cdot\nabla r^1=0,\quad \dfrac{\partial\bar{H}}{\partial r^2}\cdot\nabla r^2=0,\quad |\bar{H}\left(r^1,r^2\right)|^2=\mbox{const.}
	\end{equation}
	Note that the Alfvén double wave solution is a known solution (e.g. \cite{Zajaczkowski}) and are a stationary wave solution.
	
	Let us now derive the double Alfvén wave solution of the mhd equation (\ref{mhd_system}) obtained from the conditional symmetry method as described in Section 5. The corresponding vector field $X$, given by (\ref{Cond_symm_vector_field}), associated with the Alfvén wave vectors $\lambda_1$ and $\lambda_2$ (given by \eqref{AlfvenElements}) and satisfying the orthogonality conditions (\ref{Cond_symm_orth_cond}), takes the normalized form
	\begin{equation}
		X_1=\dfrac{\partial}{\partial t}+\mu_1\bar{h}_1^i\dfrac{\partial}{\partial x^i},\qquad X_2=\dfrac{\partial}{\partial t}+\mu_2\bar{h}_2^i\dfrac{\partial}{\partial x^i},
	\end{equation}
	where $\mu_1$ and $\mu_2$ are arbitrary functions of $u$. Under a change of variables of the form (5.11)
	\begin{align}
	\begin{split}
	\bar{t}=t,\quad \bar{x}^1=(\bar{\alpha}^1\times\bar{h}_1)_i x^i,\quad \bar{x}^2=(\bar{\alpha}^2\times\bar{h}_2)_i x^i,\\ \rho=\rho_0,\quad p=p_0,\quad \bar{v}=\dfrac{\varepsilon\bar{H}}{\sqrt{4\pi\rho_0}},\quad \varepsilon=\pm 1,\quad \rho_0,p_0\in\mathbb{R},
\end{split}	
	\end{align}
the vector fields (6.16) can be rectified to get one vector field 
\[
X=\partial/\partial\bar{t},
\]
taking $\mu_1$ and $\mu_2$ equal to zero. The m.h.d. equations (6.9) are reduced to the following overdetermined system of PDEs
\begin{align}
\begin{split}
\rho=\rho_0,\quad p=p_0,\quad |\bar{H}\left(\bar{x}^1,\bar{x}^2\right)|^2=\mbox{const.},\\ (\bar{v}\cdot\nabla)\bar{v}=\dfrac{1}{4\pi\rho}(\bar{H}\cdot\nabla)\bar{H},\quad (\bar{v}\cdot\nabla)\bar{v}=\dfrac{\varepsilon}{\sqrt{4\pi\rho}}(\bar{H}\cdot\nabla)\bar{v},\quad\nabla\cdot\bar{v}=0,\quad \bar{v}_{\bar{t}}=0,\\ (\bar{v}\cdot\nabla)\bar{H}=(\bar{H}\cdot\nabla)\bar{v},\quad (\bar{v}\cdot\nabla)\bar{H}=\dfrac{\varepsilon}{\sqrt{4\pi\rho}} (\bar{H}\cdot\nabla)\bar{H},\quad
\nabla\cdot\bar{H}=0,\quad \bar{H}_{\bar{t}}=0.
\end{split}
\end{align}
Solving (6.18), we determine a double Alfvén wave solution
\begin{equation}
\rho=\rho_0,\quad p=p_0,\quad \bar{v}=\dfrac{\varepsilon\bar{H}\left(\bar{x}^1,\bar{x}^2\right)}{\sqrt{4\pi\rho_0}},\quad |\bar{H}\left(\bar{x}^1,\bar{x}^2\right)|^2=\mbox{const.}
\end{equation}
The Riemann invariants $\bar{x}^1$ and $\bar{x}^2$ have to satisfy weaker conditions
\begin{equation}
\dfrac{\partial\bar{H}}{\partial \bar{x}^1}\cdot\nabla \bar{x}^1+ \dfrac{\partial\bar{H}}{\partial \bar{x}^2}\cdot\nabla \bar{x}^2=0,\quad |\bar{H}\left(\bar{x}^1,\bar{x}^2\right)|^2=\mbox{const.}
\end{equation}
	than the ones found by the method of characteristics (\ref{mhd_characteristics_method}). This fact is a consequence of the fact that the conditional symmetry method does not require the compatibility conditions (\ref{Gen_met_1st_order_PDEs}). Hence, a larger class of double Alfvén wave solutions of the m.h.d. equations (6.9) is obtained through the use of this method.

 \section{Final remarks} 
 
 In this paper, we have established a relation between two approaches for the construction of Riemann $k$-wave solutions of first-order qusilinear hyperbolic systems, namely the conditional symmetry method and the generalized method of characteristics. In the papers \cite{Burnat,PeradzynskiNonlin,GrundlandInvolutivity}, the involutivity conditions for the existence of k-wave solutions obtained from the generalized method of characteristics were formulated. This method required the augmentation of the initial system of PDEs (2.1) with particular differential constraints for which the matrix of derivatives of $u(x)$ is equal to a linear combination of simple elements $\xi^s\gamma_s\otimes\lambda^s$, given by (3.1), where $\xi^s$ are scalar functions of $x$. In the present paper, we derive the criterion for the existence of k-wave solutions of (2.1) through the conditional symmetry method. We supplement the initial system (2.1) with more general differential constraints for which the matrix of derivatives of $u(x)$ has the form (5.3), where simple elements are multiplied by the non-singular $k\times k$ matrix $\phi$ of $x$. A comparison with (3.1) suggests that we cannot generally expect a reduction analogue to (5.3) because $\phi$ is no longer a scalar given by (5.4). Note that if $k\geq 2$ then (5.3) is more involved than its analogue (3.1) due to the presence of the matrix $\phi$. As a consequence of direct computations, we have shown that the conditional symmetry method possesses an advantage over the generalized method of characteristics, mainly due to relaxing one of the involutivity conditions (3.4) and (3.5), namely the second-order differential condition (3.15), on the wave functions $\lambda^s$. Hence, the conditional symmetry method in the version presented in this paper can provide larger classes of k-wave solutions of a quasilinear systems of PDEs of form (2.1) written in terms of Riemann invariants through a simpler procedure.
	
	In the past, three decades ago, an interesting and significant development took place, namely the study of the Poisson structure in the connection with hydrodynamic-type systems of the form (\ref{rk1_sol_syst}). This approach was first formulated by Tsarev \cite{Tsarev}, Dubrovin and Novikov \cite{DubrovinNovikovHamil,DubrovinNovikovHydro}. Next, several non-local generalizations of the Poisson bracket with more independent variables were formulated (see e.g. \cite{Tsarev,MokhovFerapontov,Ferapontov,Mokhov}). Several new classes of multiple Riemann wave solutions were constructed. However, the problem of interrelations between the proposed two approaches and the Hamiltonian formalism still remains open. It is worth noticing that these methods seem to be complementary. The method of Hamiltonian formalism for hydrodynamic-type systems is mainly devoted to the form of quasilinear systems of PDEs, which admit a Poisson structure for the problem of $k$-waves. On the other hand, the generalized method of characteristics and the proposed version of the conditional symmetry method are mainly focused on the construction of specific classes of $k$-wave solutions for these types of systems. Such an analysis could allow us to interpret the connection between these approaches and acquire a better understanding of the physical relevance of $k$-wave superpositions in hydrodynamic-type systems expressible in terms of Riemann invariants.

\subsection*{Acknowledgements}

AMG has been supported by a research grant from NSERC of Canada and would also like to thank the Centre de Recherches Mathématiques Université de Montréal for its warm hospitality.

\label{lastpage}
\end{document}